\documentclass[10pt,a4size,reqno]{article}

\usepackage{mydef}

\begin{document}
\title{Reduced Basis Method for the Convected Helmholtz Equation}
\author{Myoungnyoun Kim$^1$, Imbo Sim$^2$\\
\small $^{1,2}$National Institute for Mathematical Sciences, Republic of Korea\\
\small $^2$Corresponding author, Email: \texttt{imbosim@nims.re.kr}
}

\date{
}
\maketitle
\begin{abstract}
We present a reduced basis approach to solve the convected Helmholtz equation with several physical parameters. Physical parameters characterize the aeroacoustic wave propagation in terms of the wave and Mach numbers. We compute solutions for various combinations of parameters and spend a lot of time to figure out the desired set of parameters. The reduced basis method saves the computational effort by using the Galerkin projection, a posteriori error estimator, and greedy algorithm. Here, we propose an efficient a posteriori error estimator based on the primal norm.  Numerical experiments demonstrate the good performance and effectivity of the proposed error estimator.
%
\end{abstract}

\section{Introduction}

Many applications such as estimation of radar cross section \cite{Chen2012},
heat transfer phenomena with high P\'eclet number \cite{Pacciarini2014},
 propagation of wave acoustics, and so on
in physics and engineering,
 are described by  partial differential equations (PDEs) with
 proper boundary conditions,
\begin{subeqnarray}
\label{eq:pde}
\calL u &=& f,\quad\mbox{in }\Om,
\\
\calB u &=& g,\quad\mbox{on }\p\Om,
\end{subeqnarray}
where
$\calL$ and $\calB$ are operators for functions on $\Om$ and $\p\Om$,
respectively, and
$\Om$ is the domain of the problem with a boundary $\p\Om$.

To reflect the physical and geometrical changes and evaluate their effect on the result,
we introduce
input parameters and outputs of interest which are just parameters and
functional values of solutions.
 Input parameters are divided into physical
 and domain parameters.
 The change of physical parameters such as
 density, porosity, frequency, absorption coefficient, flow rate, etc.\
depending on the problem,
corresponds to
 the change of the operators from
 $\calL$, $\calB$, $f$, $g$ to
 $\calL_\mu$, $\calB_\mu$, $f_\mu$, $g_\mu$,
 where the subscript $\mu$ denotes the parameter.
 The deformation of the geometrical configuration
  caused by varying of domain
 parameters \cite{Jaeggli2014} is also studied by
 the geometric parametric variation \cite{rozza2008}
 of the domain and its boundary  denoted by
  $\Om_\mu$ and $\p\Om_\mu$.
The dependence on parameters
leads us into
a parametrized partial differential equation (\PPDE)
from  \eqref{eq:pde}
\begin{subeqnarray}
\label{eq:p2de}
\calL_\mu u_\mu &=& f_\mu,\quad\mbox{in }\Om_\mu,
\\
\calB_\mu u_\mu &=& g_\mu,\quad\mbox{on }\p\Om_\mu,
\end{subeqnarray}
and a functional value
$s_\mu = l(u_\mu)$,
where $l$ is a functional of interest, and
 $u_\mu$ denotes the solution depending on the parameter.
The output can be statistical when the input is stochastic
as treated in
\cite{BBLMNP:2010,BBMNT:2009,Hassdonk2013}.

Among many aspects to view the \PPDE,
there are two main contexts, so called the real time and
 many query contexts \cite{PR:2007,rozza2008},
to be considered crucial at least in computational engineering.
The former is found in parameter estimation
or control problem, interpreted as ``deployed'' or  ``in the field''
or ``embedded.'' That is, the parameter must be estimated
rapidly ``on site''.
Meanwhile, the latter is pursued in design optimization
or multi-scale simulation.
The state equations should be solved for many parameters
\cite{Manzoni2012,Negri2013,Negri2015} in the optimization problem,
and
many calculations of small scale problems are required
to predict the macro scale properties
in the multi-scale simulation.
Following these contexts, the \PPDE\
should be solved rapidly
without severe loss of reliability, that is, while keeping the
almost same order of approximation, the evaluation must
be done as soon as possible.

According to \cite{Quarteroni2011,rozza2008},
we can regard the set of solutions generated by
parameters in a parameter domain as a smooth low-dimensional
manifold in the approximate space.
The reduced basis method (RBM) is based on the
low order approximation of the manifold owing to
the
low dimensionality of the solution manifold. Under some
sufficient assumptions, the computational task of the \PPDE\
is decomposed into the off-line and on-line stages.
In the parameter independent off-line
stage, a heavy computation is done to generate a reduced basis.
In the on-line stage, the computation for new parameter is
performed by the Galerkin projection into the reduced basis space.
The marginal number of computations gets important
since it says about the
minimum number of computations by the usual method
 to exceed the total
cost of the RBM
due to the off-line stage.
Because of the invention of a posteriori error estimators, rigorous error bounds
\cite{Chen2010}
for outputs
of interest, and effective sampling strategies,
the RBM evaluates the reliable output for many combinations of parameters
in high dimensional parameter space rapidly, which means that
the marginal number of computation gets smaller.
The reliability of the result by the RBM is
guaranteed by theoretical results in
\cite{binev:convergence,buffa,maday}.

One sufficient assumption to ensure the decomposition of the computational task
 is the affine dependence
\cite{Quarteroni2011,rozza2008}
of the \PPDE, i.e., the related forms
are expressed by
the linear combination of
parameter independent forms with
parameter dependent coefficients.
Under this assumption,
the
error
bound has many terms depending only on the dimension of the approximate space
which are independent of the parameters, and
computed during
the off-line stage.
This is good point,
but
there are two bottlenecks in the computational point of view.
Firstly,
the error bound formula is very sensitive to round-off errors, which
may show a little bit large discrepancy between
the a posteriori error bound and the on-line
efficient formula.
Secondly, the RBM is intrusive, which means that
computation of the
solutions requires intervening the matrices assembly
routines of the code.
To remove the intervention,
one can use
the empirical
interpolation method \cite{barrault,Drohmann2012,Wirtz2014}
which separates the parameter and the space variable of the
affine coefficients.

In this paper, we describe
the propagation of acoustic waves in a subsonic uniform flow
 by the time harmonic linearized Euler equation
and transform it to
a convected Helmholtz equation
 for the pressure field in \secref{sec:conv_helm}.
 The problem of the convected Helmholtz equation is
well posed  when appropriate boundary conditions are imposed,
see \cite{becache} for details.
We present a RBM for solving the convected Helmholtz equation
with
these two physical parameters.
Physical parameters
are
the Mach
and
wave numbers, which
are
 the ratio
 of the mean flow velocity and
 frequency
 to the sonic speed in the flow,
 respectively.
The outline of
the RBM is presented
with
the greedy algorithm in the pseudo code style
in \secref{sec:rbm}.
 We present numerical simulations by varying
 physical parameters
 in \secref{sec:nr}.
Finally,
several conclusions and future works are addressed
on the convected Helmholtz equation
with several parameters.

\section{Convected Helmholtz Equations}
\label{sec:conv_helm}
\subsection{Bounded domain}
We consider  compressible flows induced in a uniform subsonic flow in the direction of $x_1$
with Mach number $0 \leq \machnumber < 1$
for $(x_1,x_2)\in \mathbb{R}^2$.
Assume that  perturbations in the density $\rho$, the pressure $p$ and
all components of the velocity vector $\mathbf{u} := (u_1,\, u_2)$ are small,
and all sources and initial disturbances  bounded to the rectangular domain
\begin{equation*}
\Om= [-a_1, a_1]\times[-a_2, a_2] , \quad  a_1, a_2 > 0.
\end{equation*}
After nondimensionalizing appropriately, the flow is governed by
the linearized Euler equation
\begin{subeqnarray}\label{eq:lee}
D_t \mathbf{u}     + \nabla p & =& 0, \slabel{eq:lee1} \\
\slabel{eq:lee2}
D_t p      + \nabla \cdot \mathbf{u} & = &0, \\
\slabel{eq:lee3}
D_t \rho  + \nabla \cdot \mathbf{u} & =& 0,
\end{subeqnarray}
where $D_t = \partial_t + \machnumber \partial_{x_1}$ is
the convected derivative or the material derivative in the direction of $(M,0)$,
see \cite{hu:1996} for a detailed derivation of the equation.
Applying  $D_t$ to (\ref{eq:lee2}) and $\nabla \cdot \,$ to (\ref{eq:lee1}),
and subtracting between them
yields the convected wave equation
\begin{equation}
\label{eq:conv_wave_eq}
    \left( \partial_t  + \machnumber \, \partial_{x_1} \right)^2 p  = \Delta \, p .
\end{equation}
The Fourier transform of $(\ref{eq:conv_wave_eq})$ with respect to time
$t$
gives the  convected Helmholtz equation
\begin{equation}
\label{eq:conv_helm_eq}
(1-\machnumber^2) \frac{\partial^2 \hat{p} }{\partial x_1^2} +
\frac{\partial^2\hat{p} }{\partial x_2^2} +
2i \wavenumber \machnumber \frac{\partial \hat{p} }{\partial x_1} +
\wavenumber^2 \hat{p}  = 0.
\end{equation}
Usually, we impose a proper
boundary condition to solve
\eqref{eq:conv_helm_eq}.
For notational convenience,
$p$ is used instead of $\hat p$, then
after enforcing a general function $f$ on the right-hand side of
\eqref{eq:conv_helm_eq},
it
takes
the following divergence form
\begin{equation}
\label{eq:conv_helm_eq:div}
\Div \paren{\calM \Grad p + \calB p} +
\wavenumber^2 p = f,
\end{equation}
where
\begin{equation*}
\calM=
\brack{
\begin{array}{cc}
1-\machnumber^2 & 0
\\
0& 1
\end{array}
},
\quad
\calB
=
\brack{
\begin{array}{c}
2 i \wavenumber \machnumber
\\
0
\end{array}
}.
\end{equation*}
In one parameter problem of \eqref{eq:conv_helm_eq:div},
the wave number
$\wavenumber$ changes under a fixed Mach number $M$.
Both $M$ and $k$ varies in their domains of parameters
in two parameters problem.
In this paper, we consider
$\wavenumber$ or $\paren{\wavenumber,M}$
as a parameter $\vecmu$.
 The variational problem of the convected Helmholtz equation
 \eqref{eq:conv_helm_eq:div}
 is to find $p\in H^1(\Om)$ such that
 for given $0\le \machnumber< 1$ and $\wavenumber >0$,
 \begin{eqnarray}\label{eq:conv_helm_eq:variational}
 \quad
 -\int_\Om
	\paren{\calM \Grad p+\calB p} \cdot \Grad v\,dx
+\wavenumber^2 \int_\Om p v\,dx
=
\int_\Om f v\,dx,
	\quad
	 \mbox{for all } v\in H^1(\Om),
 \end{eqnarray}
where
 $\vecn$ is the outer normal vector.

 \subsection{Unbounded domain}
As in \cite{park2015},
we use the following notations
\begin{eqnarray*}
\Omb&=&\Set{(x_1,x_2)\in\R^2:x_-<x_1<x_+,\, -d<x_2<d}\backslash B,
\\
\Om_L^\PML&=&\Set{(x_1,x_2)\in\R^2:x_- -L<x_1\le x_-,\, -d<x_2<d}
\\
 \Om_R^\PML&=&\Set{(x_1,x_2)\in\R^2:x_+\le x_1 < x_+ + L,\, -d<x_2<d},
\end{eqnarray*}
where $B$ is an obstacle such as a circular or elliptical hole,
$x_{\pm}\in\R$
and
 $L>0$.
From \cite{park2015},
a PML formulation for the convected Helmholtz equation is
\begin{equation}
\label{eq:conv_helm_eq:pml}
\paren{1-\machnumber^2}
\paren{
\damping(x_1)
\paren{\f{\p}{\p x_1}
+
\f{i\wavenumber\machnumber}{1-\machnumber^2}}
}^2 p
+
\f{\p^2 p}{\p x_2^2}
+
\f{\wavenumber^2}{1-\machnumber^2} p = f,
\end{equation}
in $\Omt=\Om_b\cup\Om_L^\PML\cup\Om_R^\PML$.
 Here, the damping function is of the form
 \begin{eqnarray}
 \label{eq:damping}
\damping(x_1)
&=&\f{-i\frequency}{-i\frequency+\sigma(x_1)},
\\
{\sigma(x_1)}
&=&
{\sigma_0}
\paren{
(x_1-x_-)^2\chi_{(x_- -L,x_-)}(x_1) + (x_1-x_+)^2\chi_{(x_+,x_+ +L)}(x_1)},
\end{eqnarray}
where $\sigma_0$ is a parameter for the magnitude of damping
and $\chi_A(x)$ is the characteristic function on the set $A\subset\R$.
See \cite{hu:1996} for other type of the PML condition.
The divergence form of
\eqref{eq:conv_helm_eq:pml} is
\begin{equation}
\label{eq:conv_helm_eq:pml:div}
\Div\paren{\calM_\damping \Grad p +
\calB_\damping p}
+
\paren{
\wavenumber_\damping^2
-ikM \f{\p\damping}{\p x_1}
}p
=f,
\quad
\mbox{in }
\Omt,
\end{equation}
where
\begin{equation*}
\calM_\damping=
\brack{
\begin{array}{cc}
\paren{1-\machnumber^2}\damping & 0
\\
0& {\damping}^{-1}
\end{array}
},
\,
\calB_\damping
=
\brack{
\begin{array}{c}
2 i \wavenumber \machnumber\damping
\\
0
\end{array}
},
\,
\wavenumber_\damping^2
=
\f{\wavenumber^2\paren{{\damping}^{-1}-\damping \machnumber^2}}{1-\machnumber^2}.
\end{equation*}
Note that \eqref{eq:conv_helm_eq:pml:div} is the same as
\eqref{eq:conv_helm_eq:div} in the region $\Omb$ from
the definition  of the damping function.
And the variational form of \eqref{eq:conv_helm_eq:pml:div} is
 \begin{equation}\label{eq:conv_helm_eq:pml:variational}
 -\int _{\Omt}
	\paren{\calM_\damping \Grad p+\calB_\damping p} \cdot \Grad v\,dx
-ikM
 \int _{\Omt}
 \f{\p\damping}{\p x_1}  p v\,dx
+
\wavenumber_\damping^2
 \int_{\Omt}
 p v\,dx
=
\int_{\Omt}
f v\,dx,
\end{equation}
for all $v\in H^1(\Omt)$,
see \cite{park2015} for details.
The equation
\eqref{eq:conv_helm_eq:pml:variational} is also the same as
\eqref{eq:conv_helm_eq:variational} when
the support of the test function $v$ is in $\Omb$.

 \section{Reduced Basis Method}
 \label{sec:rbm}
In general,
the
RBM constructs the reduced basis using the greedy
algorithm and  precompute the parameter independent parts of matrices  at the
off-line stage. We assemble the
matrices using the
coefficients at new parameter,
solve the system and compute the output
at the on-line stage.
In the whole process,
 we restrict
the approximate space to the much smaller subspace chosen by the greedy algorithm and
discard the unnecessary modes during the calculation of the basis.
The a posteriori estimator
measures errors of approximation and
is the key to the model order reduction
\cite{barrault,binev:convergence,buffa,Chen2010,Drohmann2012,maday,Quarteroni2011,rozza2008,Wirtz2014}.
The former is given by the
property of the approximate space and chosen under the proper assumption.
The latter depends
on the reduced basis subspace.
 \subsection{Primal and Dual Problems}

Let $X$ be $H^1_0(\Om)$ with the inner product $\inner{\cdot}{\cdot}_X$
and its associated norm $\norm{\cdot}_X$.
Let $\mu$ be a parameter
selected from a certain parameter set $\pspace$.
We solve the parametrized variational form for \eqref{eq:p2de} such as
\begin{equation*}
\aform{\usol(\mu),v;\mu}=f(v;\mu),\quad \mbox{for all } v\in X
\end{equation*}
where $\aform{\cdot,\cdot;\mu}$ and $f(\cdot;\mu)$ are
  bilinear  and  linear forms depending on the parameter vector $\mu$, respectively.
We evaluate the quantity of interest $\ovar(\mu)$ as the value of a linear functional $l\in X'$
at the solution
$\usol(\mu)$
\begin{equation*}
\ovar(\mu)=l\paren{\usol(\mu);\mu}.
\end{equation*}
The
finite dimensional approximation $\usol^\calN(\mu)$ of
$\usol(\mu)$
in a smaller
function space $X^\calN\subset X$ of dimension $\calN$
satisfies
\begin{equation}\label{eq:var:approx}
\aform{\usol^\calN(\mu),v;\mu}=f(v;\mu),\quad \mbox{for all }v\in X^\calN,
\end{equation}
and
its
quantity of interest $\ovar^\calN(\mu)$ is
\begin{equation*}
\ovar^\calN(\mu)=l\paren{\usol^\calN(\mu);\mu}.
\end{equation*}
The approximate solution $\usol^\calN(\mu)$ of
\eqref{eq:var:approx}
 is  the truth approximation,
 which is accurate enough for all parameters $\mu\in\pspace$.
 To claim the accuracy, we must
 choose a very large $\calN$ and
 thus
need to solve
a large sparse matrix system of algebraic equations.

In the RBM,  we want to make a much smaller space $X_N$ than
 the approximate space
 $X^\calN$.
The space $X_N$ is called a
reduced basis,
spanned by
the linearly independent approximate solutions $\tSet{\usol^\calN(\mu_j)}_{j=1}^N$, i.e.,
$X_N=\tSpan{\tSet{\usol^\calN(\mu_j)}_{j=1}^N}$.
For the user-chosen parameter $\mu\in\pspace$ ,
the reduced basis approximation $\usol_N(\mu)\in X_N$ is obtained by the Galerkin projection,
\begin{equation}\label{eq:var:reduced}
\aform{\usol_N(\mu),v;\mu}=f(v;\mu),\quad \mbox{for all }v\in X_N,
\end{equation}
and its
quantity of interest $\ovar_N(\mu)$ is
\begin{equation*}
\ovar_N(\mu)=l\paren{\usol_N(\mu);\mu}.
\end{equation*}
Note that
the reduced basis space $X_N$ of dimension $N$  is much smaller
than the finite
approximate space  $X^\calN$ of dimension $\calN$.

To
improve the
order of convergence of output, i.e., quantity of interest $\ovar(\mu)$,
we introduce the dual problem of the primal problem \eqref{eq:var:approx}:
find $w^\calN(\mu)\in X^\calN$ such that
\begin{equation}\label{eq:var:approx:dual}
\aform{v,w^\calN(\mu);\mu}=-l\paren{v;\mu},\quad
\mbox{for all } v\in X^\calN.
\end{equation}
Its reduced basis approximation $w_N(\mu)$
 of $w^\calN(\mu)$
 is also defined by
the Galerkin projection
\begin{equation}\label{eq:var:reduced:dual}
\aform{v,w_N(\mu);\mu}=-l\paren{v;\mu},\quad
\mbox{for all } v\in X_N.
\end{equation}
Formally,
the error and residual relations of the primal problem are written as
\begin{eqnarray*}
e^\pr(\mu)
&=&\usol^\calN(\mu)-\usol_N(\mu),\\
r^\pr\paren{v;\mu}
&=&f(v;\mu)-\aform{\usol_N(\mu),v;\mu}
=\aform{\usol^\calN(\mu),v;\mu}-\aform{\usol_N(\mu),v;\mu}
\\
&=&\aform{e^\pr(\mu),v;\mu},\quad
\mbox{for all } v\in X_N,
\end{eqnarray*}
and those of the dual problem are expressed by
\begin{eqnarray*}
e^\dual(\mu)&=&w^\calN(\mu)-w_N(\mu),\\
r^\du\paren{v;\mu}&=&-l(v;\mu)-\aform{v, w_N(\mu);\mu}
=\aform{v,w^\calN(\mu);\mu}-\aform{v,w_N(\mu);\mu}
\\
&=&\aform{v,e^\dual(\mu);\mu},\quad
\mbox{for all } v\in X_N.
\end{eqnarray*}
We call $e^\pr(\mu)$, $r^\pr(\cdot;\mu)$, $e^\dual(\mu)$ and $r^\du(\cdot;\mu)$
the primal error, the primal residual, the dual error and the dual residual, respectively.
As in \cite[Section~2]{pomplun},
the dual corrected output $\ovar_N^\pd(\mu)$
is defined by
\begin{equation}\label{eq:output:primal:dual}
\ovar_N^\pd(\mu) = l\paren{\usol_N(\mu);\mu}-r^\pr\paren{w_N(\mu);\mu}.
\end{equation}
Then the error
$\ovar^\calN(\mu)-\ovar_N^\pd(\mu)$ is expressed in terms of the dual residual of the primal error,
\begin{eqnarray*}
\ovar^\calN(\mu)-\ovar_N^\pd(\mu)
&=&
l\paren{\usol^\calN(\mu);\mu}-l\paren{\usol_N(\mu);\mu}+r^\pr\paren{w_N(\mu);\mu}
\nonumber
\\
&=&
l\paren{e^\pr(\mu);\mu}+\aform{e^\pr(\mu),w_N(\mu);\mu}
\nonumber
\\
&=&
-\aform{e^\pr(\mu),w^\calN(\mu);\mu}+\aform{e^\pr(\mu),w_N(\mu);\mu}
\nonumber
\\
&=&
-\aform{e^\pr(\mu),e^\dual(\mu);\mu}
\\
&=&
-r^\du\paren{e^\pr(\mu);\mu},
\end{eqnarray*}
and it is bounded by
the norms of the primal error and
the dual residual
\begin{equation*}
\Abs{\ovar^\calN(\mu)-\ovar_N^\pd(\mu)}\le
\norm{r^\du\paren{\cdot;\mu}}_{{X}'}
\norm{e^\pr(\mu)}_{X},
\end{equation*}
where the dual norm $\norm{l}_{{X}'}$ of
any linear functional $l\in X'$
is
defined
in the usual sense:
\begin{equation*}
\norm{l}_{{X}'}
=
\sup_{v\in X} \frac{l\paren{v}}{\norm{v}_{X}}.
\end{equation*}
Note that there is
improvement
in the
convergence by the solution of a dual problem, see \cite[Section~11]{rozza2008} for more details.
To treat the non-coercive problem, we may assume that
the bilinear form of the system satisfies an inf-sup condition.
The non-zero inf-sup stability constant $\beta(\mu)$ of $\aform{\cdot,\cdot;\mu}$,
\begin{eqnarray*}
&&
\beta(\mu) = \inf_{\usol\in X} \sup_{v\in X} \frac{\Abs{\aform{\usol,v;\mu}}}{ \norm{u}_X \norm{v}_X }\ne 0
\\
&\iff&
\norm{\usol}_X\le\frac1{\beta(\mu)}\sup_{v\in X} \frac{\Abs{\aform{\usol,v;\mu}}}{ \norm{v}_X },
\quad
\mbox{for all } \usol\in X,
\end{eqnarray*}
 makes it possible to bound the norm of the primal error by the dual norm of the primal residual
\begin{equation*}
\norm{e^\pr(\mu)}_{X}
\le
\frac1{\beta(\mu)}\sup_{v\in X} \frac{\Abs{\aform{e^\pr(\mu),v;\mu}}}{ \norm{v}_X }
=
\frac1{\beta(\mu)}\sup_{v\in X} \frac{\Abs{r^\pr\paren{v;\mu}}}{ \norm{v}_X }
=
\frac1{\beta(\mu)} \norm{r^\pr\paren{\cdot;\mu} }_{X'}.
\end{equation*}
We can also
bound
the primal output error
by the norms of the output and the primal residual,
\begin{eqnarray*}
\Abs{\ovar^\calN(\mu)-\ovar_N(\mu)}
&=&
\Abs{l\paren{e^\pr(\mu);\mu}}
\le
\norm{l\paren{\cdot;\mu}}_{X'}\norm{e^\pr(\mu)}_{X}
\nonumber
\\
&\le&
\frac1{\beta(\mu)} \norm{l\paren{\cdot;\mu}}_{X'}\norm{r^\pr\paren{\cdot;\mu} }_{X'},
\end{eqnarray*}
and the dual corrected output error
by the norms of the primal and dual residuals,
\begin{equation*}
\Abs{\ovar^\calN(\mu)-\ovar_N^\pd(\mu)}\le
\frac1{\beta(\mu)}
\norm{r^\du\paren{\cdot;\mu}}_{{X}'}
\norm{r^\pr\paren{\cdot;\mu} }_{X'}.
\end{equation*}
The Riesz representation $\errrep^\pr(\mu)\in X$ of the primal residual
$r^\pr\paren{\cdot;\mu}\in X'$
 such that
\begin{equation}\label{eq:riesz:primal}
\inner{\errrep^\pr(\mu)}{v}_X = r^\pr\paren{v;\mu},
\quad \mbox{for all } v\in X,
\end{equation}
is very useful
for the computation of the error estimators of the off-line and on-line stages
 in the RBM.
From the definition of the dual norm, we obtain that
\begin{equation*}
\norm{r^\pr\paren{\cdot;\mu}}_{X'}
=
\sup_{x\in X}\frac{r^\pr\paren{v;\mu}}{\norm{v}_X}
=
\sup_{x\in X}\frac{\inner{\errrep^\pr(\mu)}{v}_X }{\norm{v}_X}
=
\norm{\errrep^\pr(\mu)}_X.
\end{equation*}

\subsection{Matrices and computational costs}
Let $\tSet{\basis_j}_{j=1}^N$ be the orthonomalized reduced basis for
$X_N=\tSpan{\tSet{\usol^\calN\paren{\mu_j}}_{j=1}^N}$,
where
$\tSet{\mu_j}_{j=1}^N$
are parameters
selected from
 some sampling strategy.
Then we may
expand the solution $\usol_N(\mu)$ for the
parameter $\mu\in\pspace$
in terms of this reduced basis
\begin{equation}\label{eq:trial:galerkin}
\usol_N(\mu) = \sum_{j=1}^N \xi_{j}(\mu) \basis_j.
\end{equation}
Inserting this expansion into \eqref{eq:var:reduced}
and applying
a test function $\basis_k$ gives us
the $k$-th row of the $N$-dimensional system of equations: for $k=1,\ldots,N$,
\begin{equation}\label{eq:var:galerkin}
\sum_{j=1}^N \xi_j(\mu) \aform{\basis_j,\basis_k;\mu} = f\paren{\basis_k;\mu},
\end{equation}
with the following output
\begin{equation}\label{eq:output:galerkin}
s_N(\mu)=\sum_{j=1}^N \xi_{j}(\mu) l\paren{\basis_j;\mu}.
\end{equation}

Denote by $\Phi$ the matrix consisting of the reduced basis $\tSet{\basis_j}_{j=1}^N$ as
its column vectors,
\begin{equation*}
\Phi = \brack{\basis_1 \cdots \basis_N}
\end{equation*}
then due to
the orthonormal property of the reduced basis in $X$, it satisfies
\begin{equation*}
\inner{\Phi^*}{\Phi}_X=\id_N,
\end{equation*}
where $\Phi^*$ is the Hermitian of $\Phi$, and $\id_N$ is the identity matrix of order $N$.
Note that the inner product of $X$ is extended to the matrices of order $N$.
Using the coefficient vector $\xi(\mu)$ whose components are $\xi_j(\mu)$,
we may rewrite \eqref{eq:trial:galerkin}, \eqref{eq:var:galerkin}
and \eqref{eq:output:galerkin} as follows:
\begin{equation*}
\usol_N(\mu)=\Phi\xi(\mu),\quad
\Phi^* A(\mu) \Phi\xi(\mu)=\Phi^* F(\mu),
\quad
s_N(\mu)=L(\mu)^*\Phi\xi(\mu),
\end{equation*}
where $A(\mu)$, $F(\mu)$ and $L(\mu)$ are the matrices
representing $\aform{\cdot,\cdot;\mu}$, $f\paren{\cdot;\mu}$ and $l\paren{\cdot;\mu}$.
Here,
the matrices satisfy the following properties: for any $\usol=\Phi\xi$, $v=\Phi\eta\in X_N$,
\begin{equation}
\label{eq:matrix:form}
\paren{\Phi\eta}^* A(\mu) \paren{\Phi\xi} = \aform{\usol,v;\mu},
\quad
\paren{\Phi\xi}^* F(\mu) = f\paren{\usol;\mu},
\quad
L(\mu)^*\Phi\xi = l\paren{\usol;\mu}.
\end{equation}
We can define the Riesz representation $\errrep^\pr(\mu)$ of the primal residual in \eqref{eq:riesz:primal}
explicitly,
\begin{equation*}
\errrep^\pr(\mu) = \innerX^{-1}\brack{F(\mu)-A(\mu)\Phi\xi(\mu)},
\end{equation*}
where $\innerX$ is the matrix due to the inner product such that
\begin{equation*}
\paren{\Phi\eta}^* \innerX \paren{\Phi\xi} = \inner{\usol}{v}_X,
\quad
\mbox{for all } \usol=\Phi\xi,  v=\Phi\eta\in X_N,
\end{equation*}
and $\innerX^*=\innerX$ from the property of the inner product in $X$.
Then the square norm $\tnorm{\errrep^\pr(\mu)}_X^2$ is
\begin{eqnarray*}
\norm{\errrep^\pr(\mu)}_X^2
&=&
F(\mu)^* \innerX^{-1}F(\mu) - F(\mu)^* \innerX^{-1} A(\mu)\Phi\xi(\mu)
\nonumber
\\
&&
-\paren{\Phi\xi(\mu)}^* A(\mu)^* \innerX^{-1}F(\mu)
+ \paren{\Phi\xi(\mu)}^* A(\mu)^* \innerX^{-1} A(\mu)\Phi\xi(\mu).
\end{eqnarray*}
Similar to the primal problem,
let $\tSet{\basis_j^\du}_{j=1}^N$ be the reduced basis from solutions
$w^\calN(\mu)$ of \eqref{eq:var:approx:dual}
for  the same parameters $\tSet{\mu_j}_{j=1}^N$.
Let $\Phi_\du$ and $\xi_\dual(\mu)$ be the matrix for the reduced basis
$\tSet{\basis_j^\du}_{j=1}^N$ and the coefficient vector for $w_N(\mu)$
of \eqref{eq:var:reduced:dual}, then we can write
\begin{equation*}
w_N(\mu)=\Phi_\du\xi_\dual(\mu),
\quad
\Phi_\du^* A(\mu) \Phi_\du\xi_\dual(\mu) = - \Phi_\du^* L(\mu).
\end{equation*}
Using \eqref{eq:matrix:form}, the dual corrected output \eqref{eq:output:primal:dual} becomes
\begin{equation*}
s_N^\pd(\mu)=L(\mu)^*\Phi\xi(\mu)-\paren{\Phi_\du\xi_\dual(\mu)}^*\paren{F(\mu)-A(\mu)\Phi\xi(\mu)}.
\end{equation*}
And let $\errrep^\dual(\mu)$ be the Riesz representation of the dual residual
\begin{equation*}
\inner{v}{\errrep^\dual(\mu)}_X = r^\du\paren{v;\mu},
\quad
\mbox{for all } v\in X.
\end{equation*}
It is expressed as
\begin{equation*}
\errrep^\dual(\mu) = \innerX^{-1}\brack{-L(\mu)-A(\mu)^* \Phi_\du\xi_\dual(\mu)}.
\end{equation*}
with the square norm $\tnorm{\errrep^\dual(\mu)}_X^2$ is
\begin{eqnarray*}
\norm{\errrep^\dual(\mu)}_X^2
&=&
L(\mu)^* \innerX^{-1}L(\mu) +L(\mu)^* \innerX^{-1} A(\mu)^*\Phi_\du\xi_\dual(\mu)
\nonumber
\\
&&
+\paren{\Phi_\du\xi_\dual(\mu)}^* A(\mu) \innerX^{-1}L(\mu)
+ \paren{\Phi_\du\xi_\dual(\mu)}^* A(\mu) \innerX^{-1} A(\mu)^*\Phi_\du\xi_\dual(\mu).
\end{eqnarray*}

In the RBM, it is very crucial  to assume that all the related forms may
be expressed as the linear combinations of parameter independent forms with
parameter dependent coefficients, or they may be affine in the parameter:
\begin{equation}\label{eq:affine}
\left\{
\begin{array}{l}
\aform{\usol,v;\mu}=\sum_{m=1}^{M_a} \acoef_{a, m}(\mu)\aqform{\usol,v},
\\
\\
\phantom{X}\,
f(v;\mu) = \sum_{m=1}^{M_f} \acoef_{f,m}(\mu)f_m{(v)},
\\
\\
\phantom{X}\,
l(\usol;\mu) = \sum_{m=1}^{M_l} \acoef_{l,m}(\mu)l_m{(\usol)}.
\end{array}
\right.
\end{equation}
Here, $\aqform{\cdot,\cdot}$, $f_m{(\cdot)}$ and $l_m{(\cdot)}$
are parameter independent forms. Clearly,
$\acoef_{a, m}(\mu)$, $\acoef_{f,m}(\mu)$ and $\acoef_{l,m}(\mu)$
are parameter dependent coefficients.
This assumption enables us to realize an efficient off-line and
on-line splitting during the computational procedure.
The above is expressed in matrices
\begin{equation*}
A(\mu)=\sum_{m=1}^{M_a} \acoef_{a, m}(\mu) A_m,
\quad
F(\mu)=\sum_{m=1}^{M_f} \acoef_{f,m}(\mu) F_m,
\quad
L(\mu)=\sum_{m=1}^{M_l} \acoef_{l,m}(\mu) L_m,
\end{equation*}
where $A_m$,  $F_m$ and $L_m$ are the matrices
representing $\aqform{\cdot,\cdot}$, $f_m\paren{\cdot}$ and $l_m\paren{\cdot}$.
When the related
forms are affine as in \eqref{eq:affine},
the approximate system \eqref{eq:var:galerkin}
becomes
\begin{equation*}
\sum_{m=1}^{M_a}\sum_{j=1}^N
\acoef_{a, m}(\mu) \xi_j(\mu) \aqform{\basis_j,\basis_k}
 =
 \sum_{m=1}^{M_f} \acoef_{f,m}(\mu) f_m\paren{\basis_k}.
\end{equation*}
Let
$\frep_m$ and $\arep_{m,j}$
be the Riesz representations
of $f_m\paren{\cdot}$ and $\aqform{\basis_j,\cdot}$ such that
\begin{equation}\label{eq:riesz:f:a}
\inner{\frep_m}{v}_X=f_m(v),\quad
\inner{\arep_{m,j}}{v}_X=\aqform{\basis_j,v},\quad
\mbox{for all } v\in X.
\end{equation}
Then we have the following representation of $\errrep^\pr(\mu)$,
\begin{equation*}
\errrep^\pr(\mu)
=
\sum_{m=1}^{M_f} \acoef_{f,m}(\mu) \frep_m
-
\sum_{m=1}^{M_a}\sum_{j=1}^N
\acoef_{a, m}(\mu) \xi_j(\mu) \arep_{m,j},
\end{equation*}
and its norm may be expressed as
\begin{eqnarray*}
\norm{\errrep^\pr(\mu) }_X^2
&=&
\sum_{m,n=1}^{M_f} \acoef_{f,m}(\mu) \bar\acoef_{f,n}(\mu)
\inner{\frep_m}{\frep_n}_X
\nonumber
\\
&&-
\sum_{m=1}^{M_f}\sum_{n=1}^{M_a}\sum_{j=1}^N
\acoef_{f,m}(\mu) \bar\acoef_{a, n}(\mu) \bar\xi_j(\mu)
\inner{\frep_m}{\arep_{n,j}}_X
\nonumber
\\
&&-
\sum_{n=1}^{M_f}\sum_{m=1}^{M_a}\sum_{j=1}^N
\bar\acoef_{f,n}(\mu) \acoef_{a, m}(\mu)  \xi_j(\mu)
\inner{\arep_{m,j}}{\frep_n}_X
\nonumber
\\
&&
+
\sum_{m,n=1}^{M_a}\sum_{j,k=1}^N
 \acoef_{a, m}(\mu) \bar\acoef_{a, n}(\mu) \xi_j(\mu) \bar\xi_k(\mu)
\inner{\arep_{m,j}}{\arep_{n,j}}_X,
\end{eqnarray*}
which is independent of $\calN$ after off-line computations of
the $\calN$ dependent quantities $\frep_m$ and $\arep_{m,j}$
with the inner products $\tinner{\frep_m}{\frep_n}_X$,
$\tinner{\frep_m}{\arep_{n,j}}_X$,
$\tinner{\arep_{m,j}}{\frep_n}_X$ and
$\tinner{\arep_{m,j}}{\arep_{n,j}}_X$.
The number of operations
to evaluate $\norm{\errrep^\pr(\mu) }_X$, or the computational cost
$\cost\paren{\norm{\errrep^\pr(\mu) }_X}$ is
\begin{eqnarray*}
\cost\paren{\norm{\errrep^\pr(\mu) }_X}
&=&
3 M_f^2+8 M_f M_a N + 5 M_a^2 N^2,
\end{eqnarray*}
where the operational costs of addition, subtraction, multiplication
and square root are assumed to be of the same order.
The coefficients $\xi_j(\mu)$ of $\usol_N(\mu)$ are obtained after solving the
reduced system
\eqref{eq:var:reduced} of dimension $N$, whose cost is denoted by
 $\cost_N$. In many cases, $\cost_N$ may not be of order $N^2$ due to
 the lack of the sparsity of the reduced system \eqref{eq:var:reduced}.

At the off-line stage,
we
solve
approximate solutions satisfying \eqref{eq:var:approx}
to form the reduced basis and orthonomalize them.
Using the reduced basis $\tSet{\basis_j}_{j=1}^N$,
we need to compute the Riesz representations of
$M_f+M_a N$ forms and
the inner products of
$M_f^2+2 M_f M_a N + M_a^2 N^2$ pairs.
Thus the computational cost $\cost_\off$ at the off-line stage is
\begin{equation*}
\cost_\off = N \cost_\calN + \cost_\res +
\paren{M_f+M_a N}\cost_\riesz+
\paren{M_f^2+2 M_f M_a N + M_a^2 N^2}(2\calN-1),
\end{equation*}
where $\cost_\calN$, $\cost_\res$,  and $\cost_\riesz$
are
the computational costs to solve the system \eqref{eq:var:approx}
of dimension $\calN$,
orthonomalize $X_N$ including a posteriori error estimators,
and
 compute the Riesz representation in \eqref{eq:riesz:f:a},
 respectively.
When the system \eqref{eq:var:approx} is sparse,
$\cost_\calN$ is of order $\calN^2$.
If the Riesz representation is bounded in $X_N$, then
$\cost_\riesz$ is of order $N$.

\subsection{Error estimator and greedy algorithm}
Examining the bounding formula of the primal error,
we may define the following error estimator and its effectivity:
\begin{eqnarray*}
\errest(\mu) = \frac1{\beta(\mu)} \norm{\errrep^\pr\paren{\cdot;\mu} }_{X},
&&
\erreff(\mu) = \frac{\errest(\mu)}{\norm{\usol^\calN(\mu)-\usol_N(\mu)}_X},
\end{eqnarray*}
where the effectivity $\erreff(\mu)$ quantifies the performance of
the error estimator $\errest(\mu)$ for the reduced basis solution.
The stability constant is assumed to be constant during the calculations,
which causes slight loss of effectivity but still works.
We
judge the current reduced basis
approximation is sufficiently accurate
if
all values of the selected error estimator are smaller than
the given tolerance.

In \algref{alg:greedy} (Greedy Algorithm)
\cite{buffa,Chen2010,Drohmann2012,Quarteroni2011,Wirtz2014},
it starts from the selection of the training sample set $\pspace_\train$
from the parameter space $\pspace$,
the tolerance $\eps$ of the error estimator
and
the maximum dimension $N_{\max}$ of the reduced basis
at \linesref{alg:greedy:train}{alg:greedy:Nmax}.
For the
randomly chosen
parameter $\mu$
 in
 $\pspace_\train$,
 we compute
the
solution of \eqref{eq:var:approx}
and normalize it with respect to the inner
product $(\cdot,\cdot)_X$ of $X$ at \linesref{alg:greedy:firstmu}{alg:greedy:firstO}.
Then we search
the next parameter maximizing the error estimator
among parameters in
 $\pspace_\train$ at \lineref{alg:greedy:secondmu}.
If the error estimator for the new parameter is smaller than the tolerance or
the dimension of the basis system is over the maximum dimension, then the
process stops at \lineref{alg:greedy:stop}. Otherwise,
we compute
the solution,
orthonormalize
the new
basis
including the previous ones,
construct
the error estimator,
find
the next
parameter maximizing the error estimator,
evaluate
the error estimator for the new
parameter,
and examine the result to decide whether to stop the process at
\linesref{alg:greedy:nextshot}{alg:greedy:nextmu}
and \lineref{alg:greedy:stop}.
Using the final reduced basis, we can compute the solution and
the residuum for new parameter at
\linesref{alg:greedy:onlinemu}{alg:greedy:onlineR}.
\algref{alg:greedy}
shows these procedure for the problem \eqref{eq:conv_helm_eq} in
 the pseudo code style.

\begin{algorithm}[H]
\caption{Greedy Algorithm}\label{alg:greedy}
\begin{algorithmic}[1]
\Procedure{Initialization:}{}
\State Construct the training sample set $\pspace_{\train}$ \label{alg:greedy:train}
\State Specify a tolerance $\eps$ as stopping criteria \label{alg:greedy:tol}
\State Choose the maximal dimension $N_{\max}$ of the reduced basis space \label{alg:greedy:Nmax}
\EndProcedure
%
\Procedure{Offline procedure:}{}
\State Choose the first parameters $\vecmu_1$ randomly \label{alg:greedy:firstmu}
\State Compute the snapshot $\vecp(\vecmu_1)$ of  \eqref{eq:conv_helm_eq} \label{alg:greedy:firstshot}
\State Set $X_1={\Set{\vecp(\vecmu_1)}}$ \label{alg:greedy:firstX}
\State Orthonormalize $X_1$ \label{alg:greedy:firstO}
\State Construct the residuum $\Delta_1$ \label{alg:greedy:firstR}
\State $\displaystyle \vecmu_2=\arg \max_{\vecmu\in \pspace_{\train}} \Delta_1(\vecmu)$
		\label{alg:greedy:secondmu}
\State $i=2$
 \While{$\Delta_{i-1}(\vecmu_i) \ge \eps$ and $i\le N_{\max}$}
 \label{alg:greedy:stop}
\State Compute the snapshot $\vecp(\vecmu_i)$ of  \eqref{eq:conv_helm_eq} \label{alg:greedy:nextshot}
\State Set $X_{i}=X_{i-1} \cup {\Set{\vecp(\vecmu_{i})}}$
\State Orthonormalize $X_{i}$ \label{alg:greedy:nextO}
\State Construct the residuum $\Delta_i$ \label{alg:greedy:nextR}
\State $\displaystyle \vecmu_{i+1}=\arg \max_{\vecmu\in \pspace_{\train}} \Delta_i(\vecmu)$
\label{alg:greedy:nextmu}
\State $i=i+1$	
\EndWhile
\EndProcedure
\Procedure{Online procedure:}{}
\State Choose new parameter $\vecmu$\label{alg:greedy:onlinemu}
\State Compute the coefficients of the reduced system
\State Determine the solution using the reduced basis
\State Compute the residuum of the solution\label{alg:greedy:onlineR}
\EndProcedure
\end{algorithmic}
\end{algorithm}

\section{Numerical Results}
\label{sec:nr}
The fundamental solution
of the convected Helmholtz equation
generated from the point source at the origin
 is
 \begin{equation*}
\Phi(x_1,x_2)
=
\f{i}{4\sqrt{1-\machnumber^2}}
H^{(1)}_{0}\paren{\f{\wavenumber\sqrt{x_1^2+(1-\machnumber^2)x_2^2}}
{1-\machnumber^2}}
\exp\paren{-i\f{\wavenumber \machnumber x_1}{{1-\machnumber^2}}},
\end{equation*}
where $H^{(1)}_{0}(z)$ is a Hankel function.

Let $K$ be a triangular element consisting of three vertices $\paren{x_1,y_1}$,
$\paren{x_2,y_2}$ and $\paren{x_3,y_3}$.
Obviously, it belongs to
the triangulation $\calT$ of $\Om$ and $\Abs{K}$ denotes its area.
The affine transformation $T^K$ from the reference triangle $\refK$ to the triangle $K$ is
defined by
\begin{equation*}
T^K:\refx\mapsto T^K(\refx) = \calB^K \refx + \vecb^K,
\end{equation*}
where
\begin{equation*}
\calB^K=
	\begin{bmatrix}
	x_2-x_1 & x_3-x_1\\
	y_2-y_1 & y_3-y_1
	\end{bmatrix},
\quad
\vecb^K=
	\begin{bmatrix}
	x_1\\
	y_1
	\end{bmatrix}.
\end{equation*}
We use the P1 conforming finite element basis function.
Then from \eqref{eq:conv_helm_eq:variational},
 the local system at the element $K$ satisfies, for a local solution vector
 $\calP^K$,
\begin{equation*}
\paren{-\calS^K+\calM^K+\calC^K}\calP^K=\calF^K,
\end{equation*}
where the stiffness $\calS^K$,  mass  $\calM^K$
and
 convection  $\calC^K$
  matrices
  are
  \begin{equation*}
\calS^K=\sum_{l=1}^3\alpha_l^K(M)\calS_l,
\quad
\calM^K=\wavenumber^2\f{\Abs{K}}{12} \calM,
\quad
\calC^K=-i \wavenumber \machnumber
\paren{\calB^K_{22}\calC_1-\calB^K_{21}\calC_2},
\end{equation*}
where the parameter independent parts are
\begin{eqnarray*}
&&
\calS_1=
	\begin{bmatrix}
	1&-1&0\\
	-1&1&0\\
	0 & 0&0
	\end{bmatrix},\quad
	\calS_2
	=
	\begin{bmatrix}
	2&-1&-1\\
	-1&0&1\\
	-1 & 1&0
	\end{bmatrix},\quad
	\calS_3
	=
	\begin{bmatrix}
	1&0&-1\\
	0 & 0&0\\
	-1&0&1
	\end{bmatrix},
	\\
	&&
\calM
=
	\begin{bmatrix}
	2&1&1\\
	1& 2&1\\
	1&1&2
	\end{bmatrix},
	\quad\quad\,\,\,
\calC_1=	
         \begin{bmatrix}
	-1&-1&-1\\
	1& 1&1\\
	0&0&0
	\end{bmatrix},\quad
	\calC_2=
	\begin{bmatrix}
	-1&-1&-1\\
	0&0&0\\
	1& 1&1
	\end{bmatrix},
\end{eqnarray*}
with parameter dependent coefficients $\alpha_l^K$ ($l=1,2,3$)
 \begin{equation*}
 \alpha_l^K(M)=
 \f{1}{4\Abs{K}}\times
\left\{
\begin{array}{ll}
{
\paren{\calB^K_{12}}^2 + \paren{\calB^K_{22}}^2
- \machnumber^2 \paren{\calB^K_{22}}^2
}
&\mbox{ if } l=1,
\\
\\
{
-\calB^K_{11}\calB^K_{12}
-\calB^K_{21}\calB^K_{22}
+\machnumber^2 \calB^K_{21} \calB^K_{22}
}
&\mbox{ if } l=2,
\\
\\
{
\paren{\calB^K_{11}}^2
\paren{\calB^K_{21}}^2
-\machnumber^2 \paren{\calB^K_{21}}^2
}
&\mbox{ if }l=3.
\\
\end{array}
\right.
\end{equation*}
Similarly, we can express the force vector $\calF^K$
after imposing appropriate boundary conditions for
the boundary integrals in \eqref{eq:conv_helm_eq:variational} such that simple Dirichlet condition.
We assemble these local systems into the global system and solve the problem for the given parameter.
For the PML case,
we can also derive a similar affine system
of the linear combination of parameter independent matrices and parameter dependent
coefficients.

The computational cost
for one realization of  uncertain parameters
in the problem
by the Galerkin method is lower than
the total computational cost including
the off-line cost $\cost_\off$ and
the on-line cost $\cost_\on$ by the RBM,
but if we want to solve the problem with
many different realizations of parameters,
the reduced basis  allows us to reduce the total
computational cost.
Let $n$ be the number of
computations due to realization of parameters.
Then
the RBM is profitable when
$n\,\cost_\galerkin \ge \cost_\off + n \,\cost_\on$,
where $\cost_\galerkin$ is the
computational
cost of the Galerkin method. In short, the profit by the RBM occurs whenever
\begin{equation}
\label{eq:marginal}
n\ge\frac{\cost_\off}{\cost_\galerkin-\cost_\on}
\end{equation}
holds.
We call
the smallest integer
$\nstar$
satisfying the inequality \eqref{eq:marginal} as the
{\em marginal number\/}
of the RBM, which indicates
the minimum number of computations to get the benefit of
the RBM in the aspect of the total computational cost.
Usually,
the computational costs are measured in seconds.

\subsection{Bounded Domain}

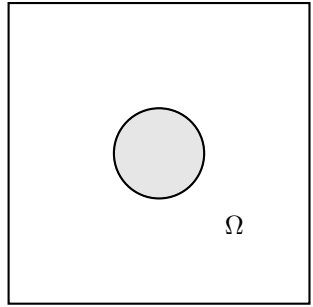
\begin{figure} [h!tbp]
\centering
	\begin{tikzpicture}[scale=2]
	\draw[thick] (-1,-1) rectangle (1,1);
	\draw[fill=gray!20] (0,0) circle [radius=0.3] ;
	\draw[thick] (0,0) circle [radius=0.3] ;
	\node[above] at (0.5,-0.6) {$\Om$};
	\end{tikzpicture}
\caption{\small Bounded domain $\Om$
excluding a circular hole.
 }
\label{fig:domain}
\end{figure}

The
bounded
domain
is
a box
$[-1, 1]\times [-1, 1]$
except
a circular hole
of radius $0.3$
and center at the origin as shown
 in \figref{fig:domain}.
We choose
between
$h=0.03$  and
$h=0.025027$
as the maximum diameter of elements in the mesh,
which is called by the mesh size.
The first choice generates
$4800$ vertices and $9250$ elements, while
the latter does  $10841$ vertices and $21157$ elements
 in the domain by \GMSH\
 \cite{GMSH:2009}.

\subsubsection{One Parameter of Wavenumber}
We use the training sample set $\pspace_\train=\Set{\wavenumber_i}_{i=1}^{N_1}$
consisting of an $N_1$ terms of an arithmetic progression sequence from
$\wavenumber_{\min}$ to $\wavenumber_{\max}$,
where $N_1$, $\wavenumber_{\min}$ and $\wavenumber_{\max}$
are the number of samples,  lower
and upper  bounds of $\wavenumber$ in \eqref{eq:conv_helm_eq}, respectively.
We
take $40$ samples ($N_1=40$) from the interval between
 $\wavenumber_{\min}=2$ and $\wavenumber_{\max}=5$.
 Numerical computations are done
 for $\machnumber=0.3$ and $\machnumber=0.4$ in the mesh of $h=0.03$.

\begin{figure} [h!tbp]
\centering
      \includegraphics[width=0.8\textwidth, trim=10mm 5mm 10mm 2mm,clip=true]
      {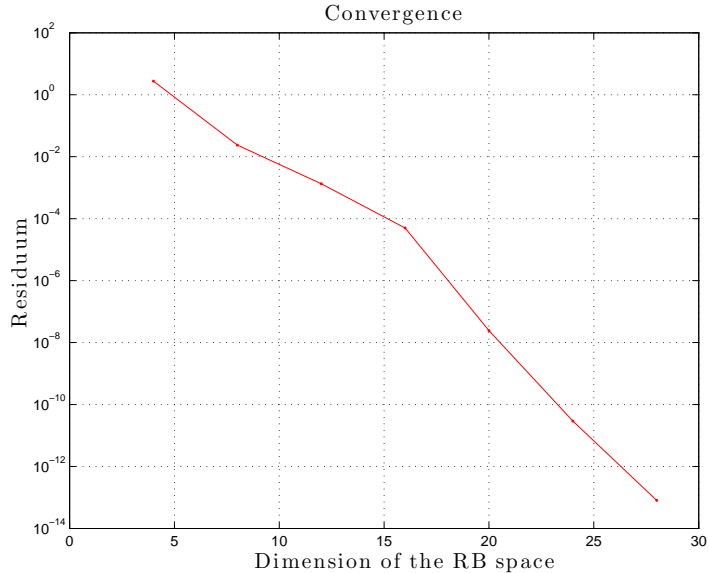}
\caption{\small Convergence in function of the dimension of the reduced basis space
 for a parameter $\wavenumber$
under
$\machnumber=0.3$.}
\label{conv1}
\end{figure}

For $\machnumber=0.3$,
\figref{conv1} comes from the residuum columns of  \tabref{tab:conv2} and
illustrates the evolution of the residuum as the dimension
$N$ of the reduced basis space  increases.
It shows
very fast decrease of residuum after the dimension
exceeds 16.
We report that
the computational costs at the on-line stage are
$\cost_\on=0.0348$ for $\machnumber=0.3$
and
$\cost_\on=0.0408$ for $\machnumber=0.4$.
Compared to the computational costs of
one  computation by the Galerkin method,
 those at the on-line stage
are $180$ and $150$ times shorter for
$\machnumber=0.3$ and $\machnumber=0.4$, respectively,
which can be calculated by taking  ratios of $\cost_\galerkin$
in \tabref{tab:conv2}
to $\cost_\on$.
The marginal number $\nstar$ increases as the dimension
of the reduced basis space does.
For instance,
the marginal numbers $\nstar$
are 876 and 871 for $\machnumber=0.3$ and $0.4$, respectively,
at the 28 dimensional reduced basis space.

\begin{figure} [h!tbp]
\centering
      \includegraphics[width=0.8\textwidth, trim=10mm 5mm 10mm 2mm,clip=true]
      {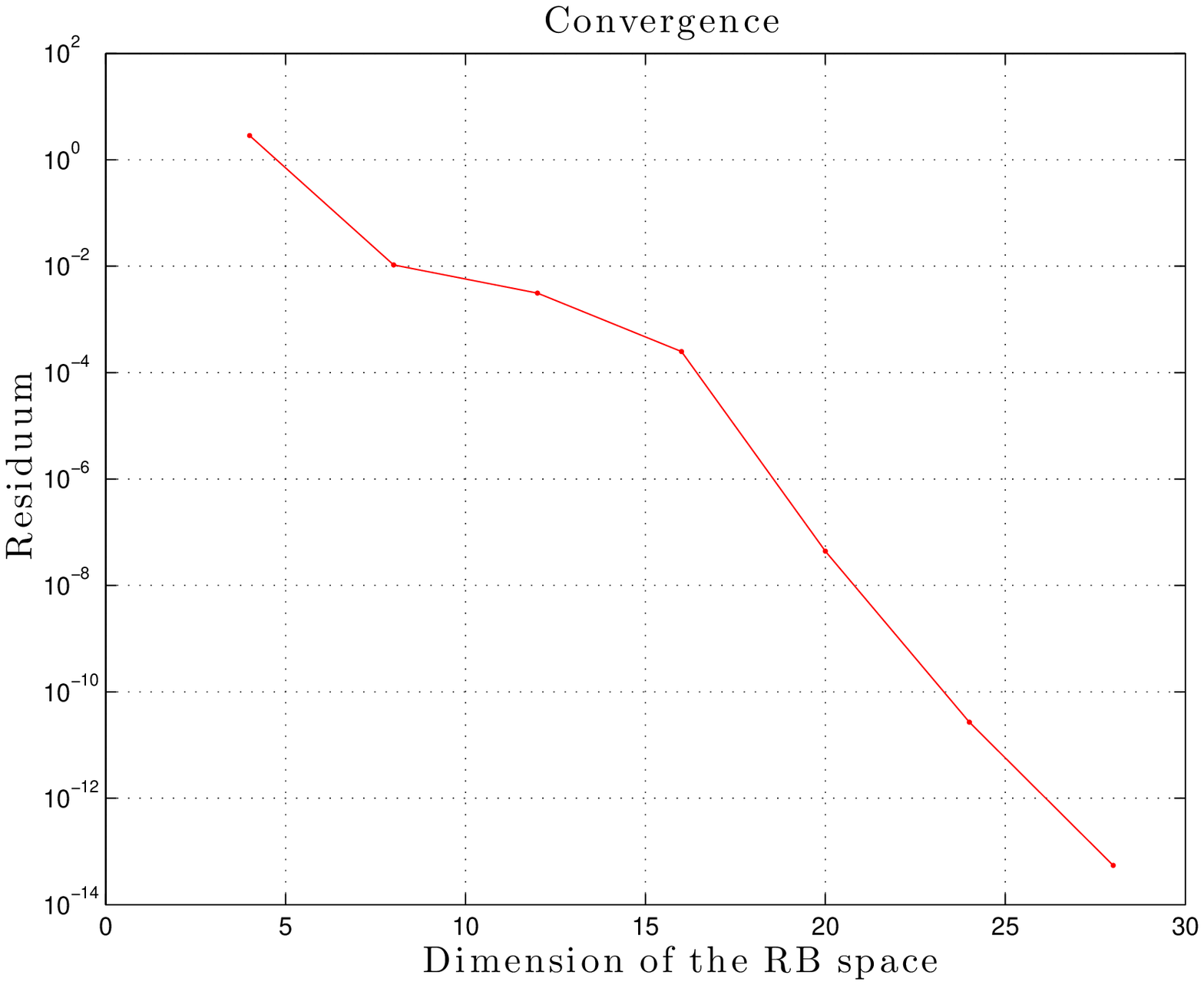}
\caption{\small Convergence in function of the dimension of the reduced basis space
 for a parameter $\wavenumber$
under
$\machnumber=0.4$.}
\label{conv1-2}
\end{figure}

\begin{table} [h!tbp]
\caption{\small Residuum, computational costs
 of the off-line stage and one Galerkin solution
 for
$\machnumber=0.3$
 and
  $\machnumber=0.4$
  according to the dimension of the reduced basis space.
  }
  \label{tab:conv2}
\centering
      \begin{tabular}{|c|l|r|c|l|r|c|}
\hline
&
\multicolumn{3}{|c|}{$\machnumber=0.3$}
&
\multicolumn{3}{|c|}{$\machnumber=0.4$}
\\
\hline
\hline
Dimension &
\multicolumn{1}{|c|}{Residuum} &
\multicolumn{1}{|c|}{$\cost_\off$} &
$\cost_\galerkin$
&
\multicolumn{1}{|c|}{Residuum} &
\multicolumn{1}{|c|}{$\cost_\off$} &
$\cost_\galerkin$ \\
\hline
4 & $2.7627$
&  789.4
& 6.37705
&    $2.8557$
&  793.9
& 6.42678
\\
\hline
8 & $2.3673\times 10^{-2}$
& 1581.1
&6.36272
&    $1.0615\times 10^{-2}$
& 1590.0
&6.41576
\\
\hline
12 & $1.3459\times 10^{-3}$
& 2370.7
&6.37737
&      $3.1121\times 10^{-3}$
& 2384.2
&6.37039
\\
\hline
16
& $4.9738\times 10^{-5}$
& 3163.7
&6.37565
& $2.4937\times 10^{-4}$
& 3181.8
&6.38359
\\
\hline
20
& $2.4432\times 10^{-8}$
& 3958.1
&6.38165
& $4.4107\times 10^{-8}$
& 3979.2
&6.42344
\\
\hline
24
& $2.9129\times 10^{-11}$
& 4750.5
&6.38063
& $2.6859\times 10^{-11}$
& 4778.2
&6.43455
\\
\hline
28
& $8.0662\times 10^{-14}$
& 5542.6
&6.37059
& $5.4560\times 10^{-14}$
& 5575.4
&6.44458
\\
\hline
\end{tabular}
\end{table}

\subsubsection{Two Parameters of Wave and Mach Numbers}
Let $\pspace_\train=\Set{\wavenumber_i}_{i=1}^{N_1}
\times\Set{\machnumber_j}_{j=1}^{N_2}$
made of
the product of
an $N_1$ terms of an arithmetic progression sequence from
$\wavenumber_{\min}$ to $\wavenumber_{\max}$,
and an $N_2$ terms of an arithmetic progression sequence from
$\machnumber_{\min}$ to $\machnumber_{\max}$.
Here $N_1$, $N_2$, $\wavenumber_{\min}$, $\wavenumber_{\max}$,
$\machnumber_{\min}$ and $\machnumber_{\max}$
are numbers of samples in wave number $\wavenumber$ and Mach number $M$ in \eqref{eq:conv_helm_eq},
lower  bounds and upper bounds of them, respectively.
We set $N_1=N_2=10$, $\wavenumber_{\max}=12$, $\wavenumber_{\min}=8$,
$\machnumber_{\max}=0.4$ and $\machnumber_{\min}=0.2$  in the mesh of
size
$h=0.025$.

The computational costs at offline $\cost_\off$,
online  $\cost_\on$ and one full Galerkin method $\cost_\galerkin$
are
23281, 0.1099 and 25.538, respectivley.
We see that
the computational benefit of the RBM occurs when the computations are more than or equal to
 the
marginal number $\nstar=916$.
We also see that
the computational cost at the on-line stage is
$230$
times shorter than that of one computation by the Galerkin method,
where the number of the speed up comes
from the ratio of the Galerkin cost $\cost_\galerkin$ to
the on-line cost.
This is very promising aspect of the RBM
such that the speed up
makes it possible
to apply the RBM to the practical problems under many and fast computational loads.

\begin{figure} [h!tbp]
\centering
\includegraphics[width=0.8\textwidth, trim=14mm 5mm 10mm 3mm,clip=true]
{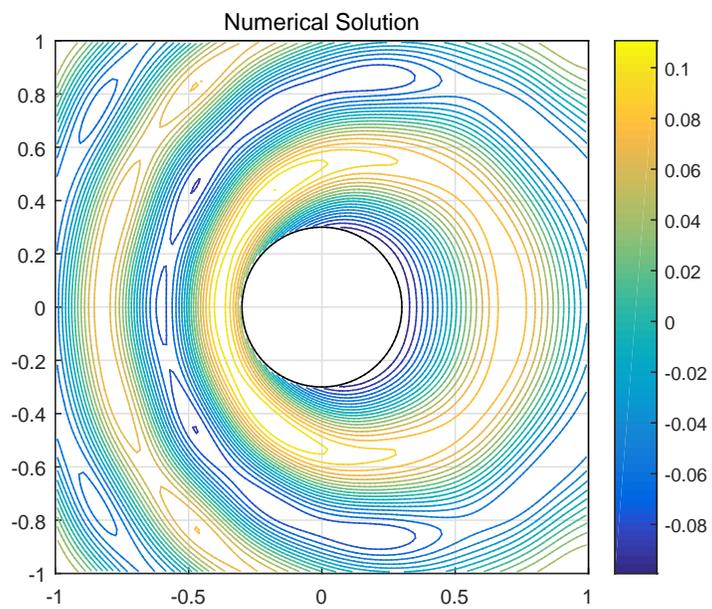}
\caption{\small
Real part of the numerical solution by the 10 dimensional RBM
chosen from 100 samples among $[8, 12]\times[0.2, 0.4]$ when
$\machnumber=0.3$ and $\wavenumber=10$.
}
\label{fig:mach3:kappa10}
\end{figure}

\figref{fig:mach3:kappa10} shows
the real part of
the
solution
 by the reduced basis of dimension $N=10$
 and the absolute error between the RBM solution and
the exact solution for fixed parameters
$\machnumber=0.3$ and $\wavenumber=10$.
The errors between the RBM solution and the exact one are
$0.0278$ in $L^\infty(\Om)$,
$0.0223$ in $L^2(\Om)$,
and
$0.0320$ in $H^1(\Om)$.

\begin{figure} [h!tbp]
\centering
\includegraphics[width=0.8\textwidth, trim=14mm 5mm 10mm 3mm,clip=true]
{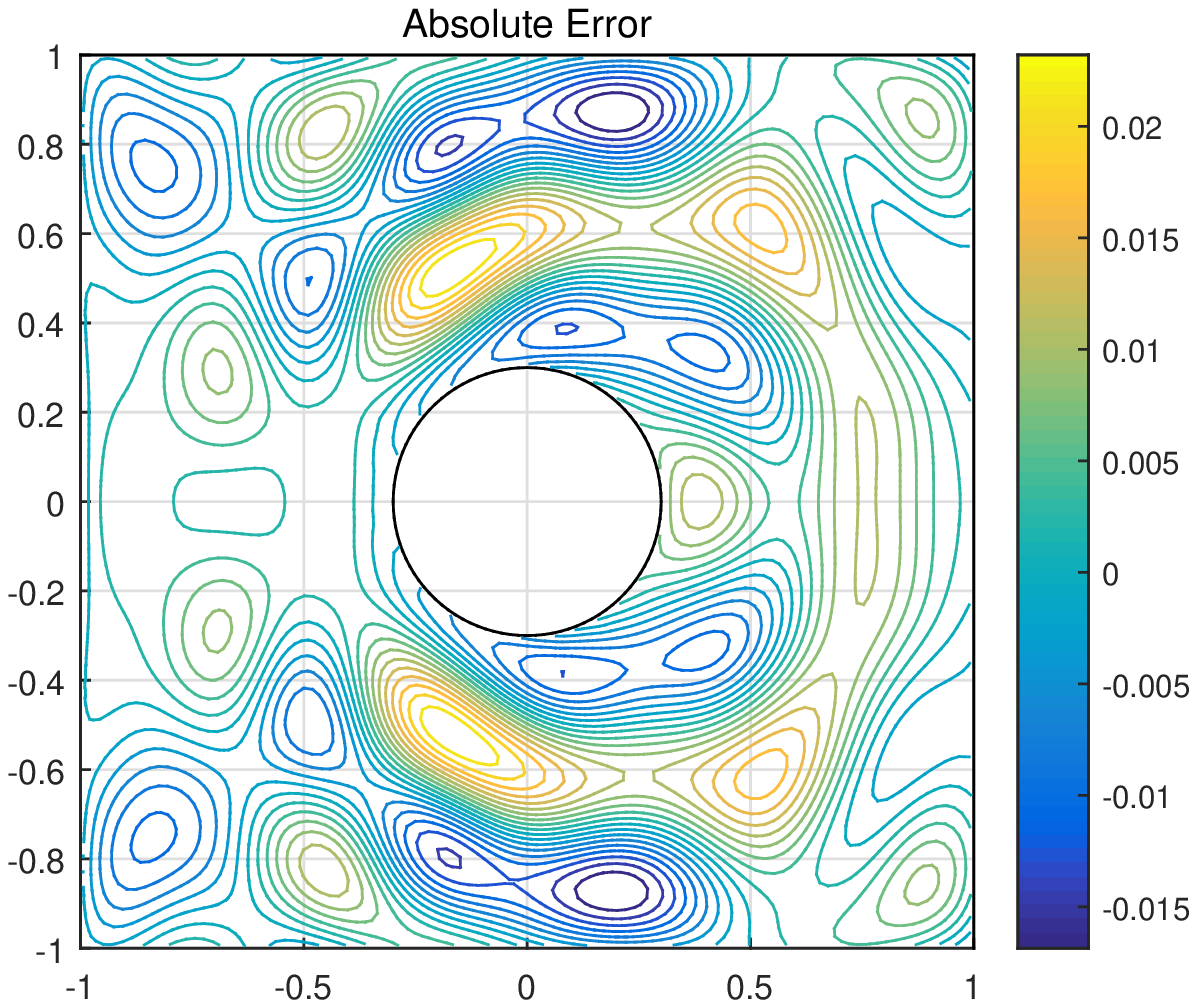}
\caption{\small
Absolute error of the numerical solution by the 10 dimensional RBM
chosen from 100 samples among $[8, 12]\times[0.2, 0.4]$ when
$\machnumber=0.3$ and $\wavenumber=10$.
}
\label{fig:mach3:kappa10:2}
\end{figure}

\subsection{Unbounded Domain}
\begin{figure} [h!tbp]
\centering
	\begin{tikzpicture}[scale=2]
	\draw[thick] (-2,-1) -- (2,-1);
	\draw[thick] (-2,1) -- (2,1);
	\draw[thick,dotted] (-1,-1) -- (-1,1);
	\draw[thick,dotted] (1,-1) -- (1,1);
	\draw[thick,dashed] (-2,-1) -- (-2,1);
	\draw[thick,dashed] (2,-1) -- (2,1);
	\draw[fill=gray!20] (0,0) ellipse (0.3 and 0.25) ;
	\draw[thick] (0,0) ellipse (0.3 and 0.25) ;
	\node[right] at (0.3,-0.6) {$\Omb$};
	\node[right] at (1.3,-0.6) {$\Om_R^\PML$};
	\node[right] at (-1.8,-0.6) {$\Om_L^\PML$};
	\end{tikzpicture}
\caption{\small Bounded domain $\Omb$ and
 PML domain $\Om_L^\PML\cup\Om_R^\PML$.
 }
\label{fig:domain3}
\end{figure}
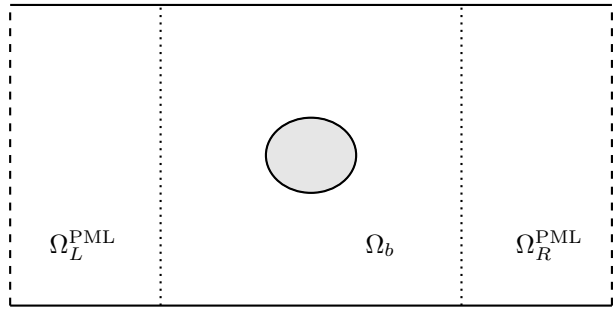

The duct
 in \figref{fig:domain3}
has
an elliptical hole whose
major and minor axes are $a=0.3$ and $b=0.25$,
and center is at the origin.
We set
$x_-=-1$, $x_+=1$, $L=1$ and $\sigma_0=15$
for
the damping function $\damping(x)$ in \eqref{eq:damping}.
We generate
meshes
for $\Omt$
of
mesh
size $h=0.0381$
by \GMSH,
which has
$16907$ nodes and $33262$ elements.
We treat the wave and Mach numbers  as  parameters.
We use 16 training samples among $[8, 12]\times[0.2, 0.4]$
and choose 10 basis from them.
The computational costs at offline $\cost_\off$,
online  $\cost_\on$ and one full Galerkin method $\cost_\galerkin$
are
4899, 0.15763 and 188.2764, respectivley.
The
marginal number is $\nstar=27$ and
the computational speed by the
on-line stage is at least 1,100 faster than that by the usual
Galerkin method.

The errors between the 10 dimensional RBM solution and the exact one are
$0.0682$ in $L^\infty(\Om)$,
$0.0248$ in $L^2(\Om)$,
and
$0.4012$ in $H^1(\Om)$.
The $H^1(\Om)$ error is higher than that for the bounded domain,
which is caused by the small number of
training samples.

\begin{figure} [h!tbp]
\centering
\includegraphics[width=0.8\textwidth, trim=14mm 5mm 10mm 3mm,clip=true]
{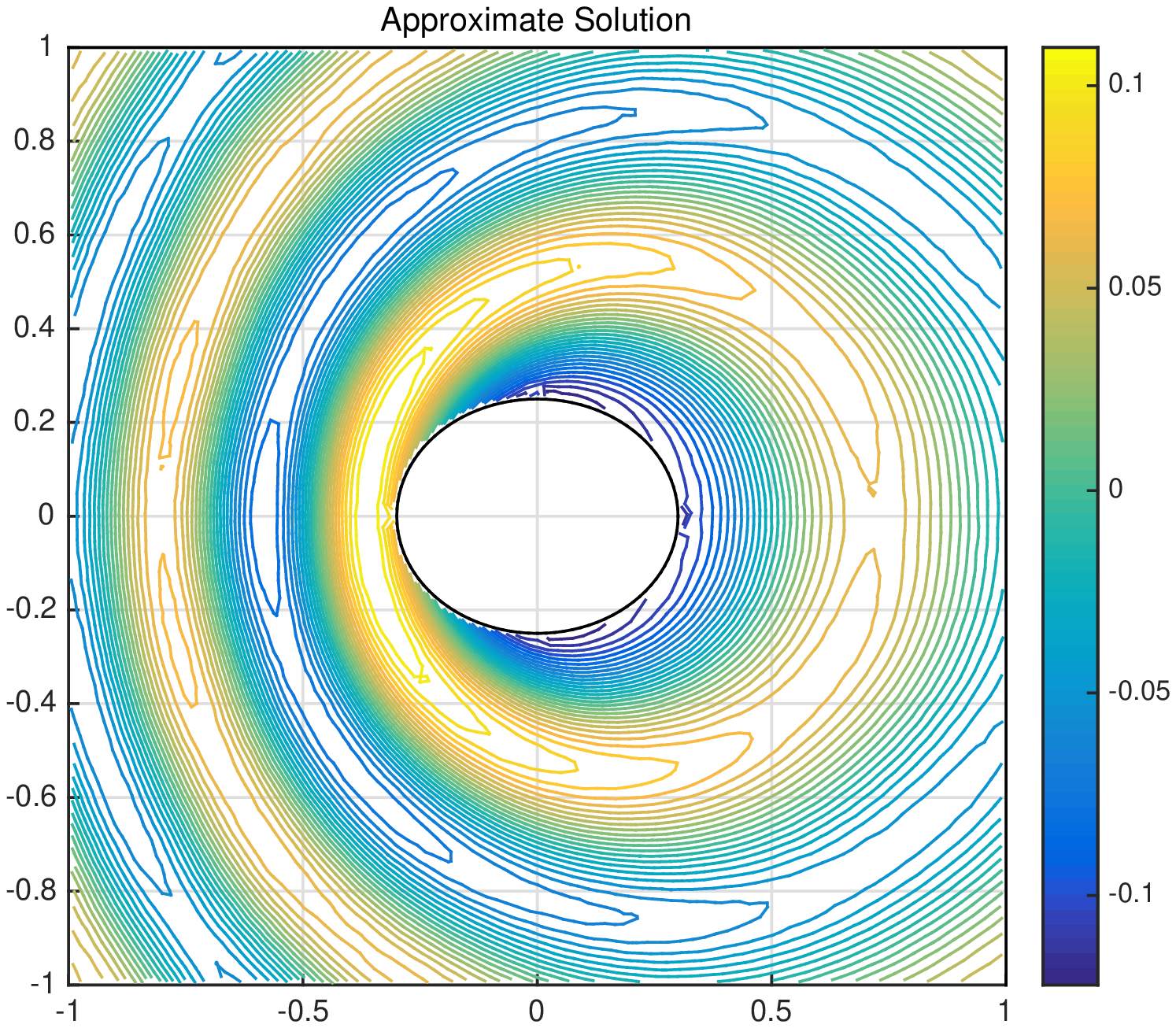}
\caption{\small
Real part of the numerical solution by the 10 dimensional RBM
selected from 16 samples among $[8, 12]\times[0.2, 0.4]$
 in $\Om_b$ when
$\machnumber=0.3$ and $\wavenumber=10$.
}
\label{fig:mach3:kappa10:ellipse}
\end{figure}

\begin{figure} [h!tbp]
\centering
\includegraphics[width=0.8\textwidth, trim=14mm 5mm 10mm 3mm,clip=true]
{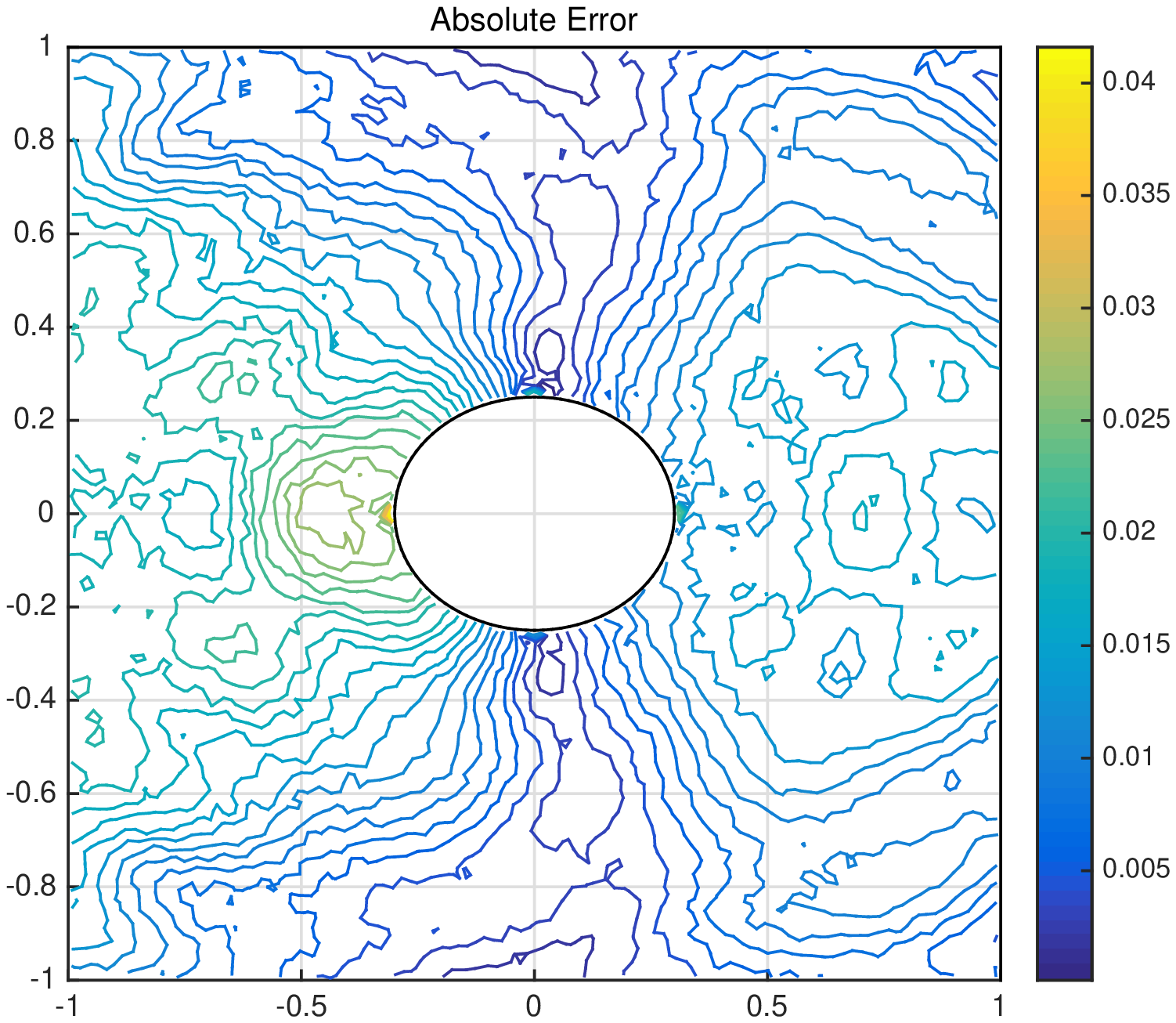}
\caption{\small
Absolute error of the numerical solution by the 10 dimensional RBM
selected from 16 samples among $[8, 12]\times[0.2, 0.4]$
 in $\Om_b$ when
$\machnumber=0.3$ and $\wavenumber=10$.
}
\label{fig:mach3:kappa10:ellipse}
\end{figure}
\section{Conclusion}
\label{sec:conclusion}
We test the RBM for
the convected Helmholtz equation.
 The physical parameters are expressed as coefficients of the equation.
 After these tests, we confirm that the RBM works well and gives
 us the benefit of fast computation at least 100 times
 than the usual computational method does.
In the implementation,
we use the error estimator based on the primal norm of the error.

\end{document}